\documentclass[11pt]{article}
\usepackage[affil-it]{authblk}

\usepackage[T1]{fontenc}
\usepackage[cp1250]{inputenc}

\usepackage[leqno]{amsmath}
\usepackage{amsfonts}
\usepackage{amssymb}
\usepackage{amsthm}
\usepackage{mathrsfs}

\usepackage{enumerate}
\usepackage{setspace}

\newcommand{\al}{\alpha}
\newcommand{\be}{\beta}
\newcommand{\ga}{\gamma}
\newcommand{\del}{\delta}
\newcommand{\Del}{\Delta}

\newcommand{\veps}{\varepsilon}
\newcommand{\lam}{\lambda}
\newcommand{\si}{\sigma}

\newcommand{\cl}[1]{\overline{#1}}

\newcommand{\bsh}{\backslash}

\newcommand{\re}{\text{Re}}

\newcommand{\power}[1]{2^{#1}} 
\newcommand{\bbC}[1][]{\mathbb{C}^{#1}}
\newcommand{\bbR}[1][]{\mathbb{R}^{#1}}
\newcommand{\bbN}[1][]{\mathbb{N}^{#1}}

\renewcommand{\H}{H} 
\newcommand{\weak}{\rightharpoonup}

\newcommand{\LR}[1][]{L^2(\mathbb{R}^{#1}\!,\bbR)}
\newcommand{\LC}[1][]{L^2(\mathbb{R}^{#1}\!,\bbC)}

\newcommand{\LL}{L} 
\newcommand{\N}{N} 
\newcommand{\dom}{\operatorname{Dom}}
\newcommand{\ran}{\operatorname{Ran}}
\renewcommand{\ker}{\operatorname{Ker}}
\newcommand{\graph}{\operatorname{Graph}}

\newcommand{\lin}{\mathcal{L}} 
\newcommand{\bound}{\mathcal{B}} 
\renewcommand{\sp}{\si}
\newcommand{\rez}{\rho} 
\newcommand{\lH}{\H_{-}}

\newcommand{\rH}{\H_{+}}
\newcommand{\lP}{P_{-}}
\newcommand{\rP}{P_{+}}
\newcommand{\lL}{\LL_{-}}
\newcommand{\rL}{\LL_{+}}
\newcommand{\lN}{\N_{-}}
\newcommand{\rN}{\N_{+}}

\newtheorem{theorem}{Theorem}

\newtheorem{lemma}[theorem]{Lemma}

\theoremstyle{definition}
\newtheorem{definition}[theorem]{Definition}
\theoremstyle{remark}
\newtheorem{remark}{Remark}

\begin{document}

\title{Solvability of semilinear equations with zero on the boundary of spectral gap and applications to nonlinear Schr\"{o}dinger equation}
\author{Przemys{\l}aw Zieli\'{n}ski
\thanks{Electronic address: \texttt{zielinski.przemek@gmail.com}}}
\affil{\em{Institute of Mathematics Polish Academy of Sciences,\\ \'{S}niadeckich 8, 00-950 Warsaw, Poland}} 
\maketitle

\begin{abstract}
We study the existence of solutions in Hilbert space $\H$ of the semilinear equation
\[
\LL u+\N(u)=h,
\]
where $\LL$ is linear self-adjoint, $N$ is a nonlinear operator and $h\in\H$. We concentrate on the case when $0$ is a right boundary point of a gap in the spectrum of $\LL$ and an element of essential spectrum. The sufficient conditions for solvability are based on monotonicity and sign assumptions on operator $\N$, and its behaviour on $\ker\LL$. We illustrate the main theorem by an application to the study of nonlinear stationary Schr\"{o}dinger equation on $\bbR[n]$.
\end{abstract}

\medskip


\section{Introduction}

Let $\H$ be a separable and complex Hilbert space with the scalar product $\langle\cdot\,,\cdot\rangle$ and norm $\|\cdot\|$. In this paper we study the solvability of the semilinear equation
\begin{equation}\label{eq:main}
\LL u+\N(u)=h,
\end{equation}
where $\LL\colon\dom(\LL)\subset\H\to\H$ is a densely defined and self-adjoint linear operator (in general unbounded), $\N\colon\H\to\H$ is nonlinear such that $N(0)=0$, and $h\in\H\bsh\{0\}$. We also make the following initial assumptions concerning the spectrum of $\LL$
\begin{enumerate}
\renewcommand{\theenumi}{$(\LL_\arabic{enumi})$}
\renewcommand{\labelenumi}{\theenumi}
\item\label{h:spectr} $0\in\sp(L)$,
\item\label{h:gap} $(-\del,0)\subset\rez(L)$\ \text{for some}\ $\del>0$,
\item\label{h:bound} $\inf\sp(\LL)\geqslant-\ga$ for some $\ga\geq 0$,
\end{enumerate}
where $\rez(\LL)\subset\bbC$ denotes the resolvent set of $\LL$
\[
\rez(\LL)=\{\lam\in\bbC:\ \LL-\lam I\ \text{is bijection from}\ \dom(\LL)\ \text{onto}\ \H\},
\]
and $\sp(\LL)=\bbC\bsh\rez(\LL)$ is the spectrum of operator $\LL$.
\begin{remark}
\begin{enumerate}[(a)]
\item Due to condition \ref{h:spectr}, the linear part $\LL$ is non-invertible. In this case we say that equation (\ref{eq:main}) is \emph{at resonance}. This is more complicated situation compared to $0\in\rez(\LL)$ since then equation (\ref{eq:main}) reduces to the fixed point problem
\[
u=\LL^{-1}(h-N(u))
\]
and a suitable topological degree theory gives a wealth of results provided $\N$ is assumed to be compact, monotone, A-proper, etc.

\item Condition \ref{h:gap}, together with \ref{h:spectr}, precise further that $0$ lies on the boundary of a spectral gap of the operator $\LL$. In particular it is possible that $0$ is an eigenvalue of $\LL$. Such a situation is typical when the linear operator arises from boundary value problem on bounded domain in $\bbR[n]$. The resolvent of such an operator, i.e., the operator-valued mapping $\rez(\LL)\ni\lam\mapsto(\LL-\lam I)^{-1}$, has compact values and hence the spectrum $\sp(\LL)$ is discrete and composed purely of eigenvalues with finite-dimensional eigenspaces, see \cite[Ch. 5.4]{BlaExnHav2008}. From the pioneering work of Landesman and Lazer \cite{LandesmanLazer1970} the literature on resonance problems for partial and ordinary BVPs in bounded domains has vastly expanded, see e.g. \cite{KarakostasTsamatos2001},\cite{Nieto1992},\cite{{Przeradzki1993}} and references therein.

\item The most important situation for us, encompassed by assumptions \ref{h:spectr} and \ref{h:gap}, arises when $0$ is an accumulation point of the spectrum $\sp(\LL)$. In this case $0$ belongs to the essential spectrum of $\LL$ (the subset of $\sp(L)$ composed of its accumulation points or eigenvalues with infinite-dimensional eigenspace). In practice the essential spectrum appears when we consider differential operators on unbounded domains in $\bbR[n]$. Among many examples we can mention: Sturm-Liouville operators 
on the half-line $[0,+\infty)$ (see \cite{Teschl2009}), the Schr\"{o}dinger operator $S=-\Del+V$ on $\bbR[n]$ (see \cite{Teschl2009},\cite{BlaExnHav2008} or below), the Dirichlet Laplacian $-\Del^D$ on periodic unbounded strips (waveguides) in $\bbR[2]$ (see e.g. \cite{Yoshitomi1998},\cite{SobolevWalthoe2002}).

The basic feature of this situation, which distinguish it from the case of operators with purely discrete spectrum, discussed in item (b), is that the values of the resolvent of $\LL$ are not compact operators any more.

\item Let us also stress out that the assumptions $\N(0)=0$ and $h\in\H\bsh\{0\}$ imply that equation (\ref{eq:main}) does not have a trivial, i.e. zero, solution. This is important since our method does not guarantee non-triviality of solutions.
\end{enumerate}
\end{remark}
\medskip

The main goal of this paper is to give the additional conditions concerning operators $\LL$ and $\N$ that will guarantee the existence of a solution of the equation (\ref{eq:main}). Since there is a lack of compactness caused by non-emptiness of the essential spectrum we will rely on monotonicity conditions and the surjectivity properties of maximal monotone operators. In this direction our main result reads as follows.

\begin{theorem}\label{thm:main_ex_monot}
Assume that a linear self-adjoint operator $\LL\colon\dom(\LL)\subset\H\to\H$ satisfies \ref{h:spectr}, \ref{h:gap}, \ref{h:bound} and that $h\in\H\bsh\{0\}$. Let $\N\colon\H\to\H$ be bounded, demicontinuous, $\N(0)=0$ and assume that the following conditions hold 
\begin{enumerate}[(i)]
\item there exists $\al>\ga/\del^2$ such that for all $u_1,u_2\in\H$
\[
\re\langle\N(u_1)-N(u_2),u_1-u_2\rangle\geqslant\al\|\N(u_1)-\N(u_2)\|^2,
\]
\item $\displaystyle\limsup_{k\to+\infty}\frac{\re\langle\N(u_k),u_k\rangle}{\|u_k\|}>\frac{\ga\|h\|}{\del^2\al-\ga}$ for each sequence $\{u_k\}_{k\in\bbN}\subset\H$ such that $\|u_k\|\to+\infty$,
\item $\displaystyle J_{\N}(u)>\frac{\ga\|h\|}{\del^2\al-\ga}+\re\langle h,u\rangle$, for all $u\in\ker\LL\bsh\{0\}$ (see definition \ref{def:recession}).
\end{enumerate}
Then the equation (\ref{eq:main}) has a solution.
\end{theorem}
\begin{remark}\label{rem:main_ex_monot}
\begin{enumerate}[(a)]
\item Note that assumption $(i)$, since the left hand side is always non-negative, implies in particular that $\N$ is monotone (definition \ref{def:monotony}(1)). Moreover, it is equivalent to the fact that the multivalued inversion of $\N$ is strongly monotone (definition \ref{def:monotony}(2)). This assumption was previously applied in the work of Brezis and Nirenberg \cite{BrezisNirenberg1978}. For example, in the case $\N$ is the gradient of a convex $C^1$ functional it is equivalent to the Lipschitz condition with constant $1/\al$, see \cite[Prop. A5]{BrezisNirenberg1978} and the remark that follows.

\item Functional $J_{\N}\colon\H\to[-\infty,+\infty]$ is the so called \emph{recession function} and is defined, roughly speaking, as (see definition \ref{def:recession})
\[
J_{\N}(u)=\liminf_{\substack{t\to\infty\\ v\weak u}}\re\langle\N(tv),v\rangle.
\]
Assumption $(iii)$ describes the behaviour of operator $\N$ on $\ker\LL$. However, in our case $0$ can be an accumulation point of the spectrum and it can have trivial eigenspace at the same time. Hence $(iii)$ cannot provide enough control on $\N$ and thus we introduce also the sign condition $(ii)$. Let us also note, that such situation never occurs when (\ref{eq:main}) is at resonance and $\LL$ has purely discrete spectrum.
\end{enumerate}
\end{remark}
Theorem \ref{thm:main_ex_monot} and the methods used in this paper are inspired by the work of Brezis and Nirenberg \cite{BrezisNirenberg1978} in which the similar recession functional was introduced and applied to semilinear problems. The authors investigated mostly the case corresponding to compact resolvent but some applications to "non-compact" problems via monotonicity methods are also provided. However, in a "non-compact" variant, only nonlinearities which are the gradients of convex functionals are treated in \cite{BrezisNirenberg1978}.
\medskip

We will give an application of theorem \ref{thm:main_ex_monot} to the nonlinear stationary Schr\"{o}dinger equation in $\bbR[n]$
\begin{equation}\label{eq:Schrodinger}
-\Del u+V(x)u+f(x,u)=h(x)
\end{equation}
where $\Del u=\sum_{i=1}^{n}u_{x_ix_i}$ is the Laplacian of $u\colon\bbR[n]\to\bbR$, $V\colon\bbR[n]\to\bbR$ is a given measurable function (potential) and $f\colon\bbR[n]\times\bbR\to\bbR$ is a Carath\'{e}odory function. We will assume that $V$ satisfies also the following condition
\begin{equation}\tag{K}\label{eq:Kato}
V\in L^{p}_{\mathrm{loc}}(\bbR[n])\ \text{and}\ V_{-}:=\min\{V,0\}\in L^{p}(\bbR[n])+L^{\infty}(\bbR[n])\quad\text{for some}\ p\in[2,\infty)\cap(n/2,\infty).
\end{equation}
This assumption ensures that linear Schr\"{o}dinger operator $S=-\Del+V$ is defined and self-adjoint in $L^2(\bbR[n])$, and its spectrum is bounded from below, see \cite[Theorem A.2.7.]{Simon1982}.

\begin{theorem}\label{thm:Schrodinger_ex}
Suppose $V$ satisfies condition (\ref{eq:Kato}) and that there is a $\del>0$ such that $(-\del,0]\cap\sp(S)=\{0\}$. Moreover assume that $f\colon\bbR[n]\times\bbR\to\bbR$ is a Carath\'{e}odory function satisfying
\begin{enumerate}[(i)]
\item for a.e. $x\in\bbR[n]$ the function $f(x,\cdot)$ is non-decreasing and $f(x,0)=0$,
\item there is a function $q\in L^2(\bbR[n])$ and constants $a,b>0$ such that for a.e. $x\in\bbR[n]$ and all $t\in\bbR$
\[
a|t|\leqslant|f(x,t)|\leqslant q(x)+b|t|,
\]
\item there is constant $\al>|\inf\sp(S)|/\del$ such that
\[
|f(x,t)-f(x,t')|\leqslant\frac{1}{\al}|t-t'| 
\]
for almost all $x\in\bbR[n]$ and all $t,t'\in\bbR$.
\end{enumerate}
Then for each $h\in L^2(\bbR[n])$ equation (\ref{eq:Schrodinger}) has a real-valued solution $u\in H^2(\bbR[n])$.
\end{theorem}
Equation (\ref{eq:Schrodinger}) with trivial solution, i.e., when $h\equiv0$, and in case when $0$ is a boundary point of spectral gap, was studied in \cite{BartschDing1999}, \cite{WillemZou2003}, \cite{YanCheDin2010} and \cite{Tang2014}. In these papers the variational methods were used to show the existence of non-trivial weak solution belonging to $H^2_{\mathrm{loc}}(\bbR[n])\cap L^s(\bbR[n])$ where $2<s<2^*$. However, the assumption that the potential $V$ and nonlinearity $f$ are periodic functions in the $x$-variable was made in all mentioned papers. In particular, when $V$ is periodic the spectrum of Schr\"{o}dinger operator $\sigma(S)$ is purely continuous (i.e., does not contain any eigenvalues), bounded from below and consists of closed disjoint intervals (so called band structure). In general, when potential is not periodic, we can have various other behaviours such as eigenvalues situated in gaps between the intervals or contained inside the intervals (the latter are called embedded eigenvalues). Theorem \ref{thm:Schrodinger_ex} is an attempt to establish solvability of equation (\ref{eq:Schrodinger}) in this more general direction.

Let us note that each function $f\colon\bbR[n]\times\bbR\to\bbR$ such that $f(x,\cdot)$ is non-decreasing and asymptotically linear, for a.e. $x\in\bbR[n]$, satisfies assumption $(ii)$ from theorem \ref{thm:Schrodinger_ex}. For an example of not asymptotically linear function $f$ consider
\[
f(x,t)=g(x,t)t,
\]
where $g\colon\bbR[n]\times\bbR\to\bbR$ is non-negative and bounded. Moreover we assume that $a\leqslant g(x,t)\leqslant1/(2\al)$ and that $\partial_tg$ exists and $|\partial_tg(x,t)|\leqslant1/(2\al|t|)$ for almost all $x\in\bbR[n]$ and $|t|\neq0$.
\medskip

\section{Preliminaries}

At first we will use the assumptions \ref{h:spectr} -- \ref{h:bound} to perform the decomposition of the Hilbert space according to the decomposition of the spectrum $\sp(\LL)$.
Let $\{E_\mu:\ \mu\in\bbR\}$ denote the spectral family of operator $\LL$. Setting
\begin{equation*}
\lH = E_{0-}\H=E_{-\del}H,\quad \rH = (I-E_{0-})\H,
\end{equation*}
where $E_{0-}=\lim_{k\to+\infty}E_{-1/k}$, we have the decomposition
\begin{equation}\label{eq:decomp}
\H=\lH\oplus\rH.
\end{equation}
Since the spaces $\lH,\rH$ are invariant with respect to operator $\LL$, the spectral theory implies the following facts.
\begin{lemma}\label{lem:decomp}
Let us denote by $\LL_{\pm}\colon\dom(\LL_{\pm})\subset\H_{\pm}\to\H_{\pm}$, where $\dom(\LL_{\pm})=\dom(\LL)\cap\H_{\pm}$, the parts of $\LL$ in $\H_{\pm}$. Then
\begin{enumerate}[(i)]
\item $\sp(\lL)=\sp(\LL)\cap[-\ga,-\del]$, i.e., the operator $\lL$ is bounded and invertible with $\|\lL^{-1}\|\leq\del^{-1}$, moreover for all $u\in\lH$ we have the estimate
\[
\langle\lL u,u\rangle\geqslant -\ga/\del^2\|\lL u\|^2;
\]
\item $\sp(\rL)=\sp(\LL)\cap[0,+\infty)$, i.e., the operator $\rL$ is non-negative.
\end{enumerate}
\end{lemma}
The proof is based on the properties of spectral integral which allow us to write
\[
\LL=\int_{-\infty}^{+\infty}\mu\,\mathrm{d}E_\mu= \int_{-\infty}^{-\del}\mu\,\mathrm{d}E_\mu+ \int_{0}^{+\infty}\mu\,\mathrm{d}E_\mu.
\]
The integrals on the right hand side define the operators $\LL_-$ and $\LL_+$ respectively.\medskip

Let us now recall some notions that are used when dealing with nonlinear maps. If we have a sequence $\{u_0\}\cup\{u_k\}_{k\in\bbN}\subset\H$ then the formula $u_k\to u_0$ will represent convergence in the norm of the space $\H$, that is
\[
\lim_{k\to\infty}\|u_k-u_0\|=0,
\]
and by $u_k\weak u_0$ we will mean the convergence of this sequence in the weak topology $\si(H,H^*)$ of $\H$, which is equivalent to the following condition
\[
\lim_{k\to\infty}\langle u_k-u_0,v\rangle=0,
\]
for every $v\in\H$.
\begin{definition}
We say that operator $\N\colon\H\to\H$ is 
\begin{enumerate}[(1)]
\item \emph{bounded} if the image of every bounded subset of $\H$ is bounded,
\item \emph{demicontinuous} if $u_k\to u_0$ implies $\N(u_k)\weak\N(u_0)$,
\item \emph{continuous} if $u_k\to u_0$ implies $\N(u_k)\to\N(u_0)$.
\end{enumerate}
\end{definition}
The demicontinuity of the mapping $\N$, although it is defined with the use of sequences, is equivalent to the continuity of $N$ from $\H$ endowed with norm topology to $\H$ with its weak topology $\si(\H,\H^\star)$. Henceforward we will denote
\[
\langle\cdot\,,\cdot\rangle_r:=\re\langle\cdot\,,\cdot\rangle.
\]

\begin{definition}\label{def:recession} 
Let $\N\colon\H\to\H$. Define the functional $J_{\N}\colon\H\to[-\infty,+\infty]$ with formula
\[
J_{\N}(u)=\inf\left\{\liminf_{k\to+\infty}\langle N(t_kv_k),v_k\rangle_r:\ t_k\to+\infty,\ \{v_k\}_{k\in\bbN}\subset\H,\ v_k\weak u\right\}.
\]
\end{definition}

\begin{remark}\label{rem:recession}
\begin{enumerate}[(a)]
\item The functional $J_{\N}$ is an example of so called sequential recession function introduced in \cite[Rem. 2.17, p. 157]{Baiocchi1988}, along with more general topological recession function \cite[Def. 2.2, p. 152]{Baiocchi1988}, to study the abstract minimization problems with non-coercive and non-convex energy functional.
\item Let $\N\colon\H\to\H$. 
Define the functional $\psi_\N\colon\H\to\bbR$ as
\[
\psi_\N(u)=\langle\N(u),u\rangle_r.
\]
Then we see that
\[
J_{\N}(u)=\liminf_{\substack{t\to\infty\\ v\weak u}} \frac{\psi_{N}(tv)}{t}.
\]
Hence, for $u\neq0$, the recession function $J_{\N}(u)$ describes the growth of $\psi_{\N}$ as we are heading to infinity in norm and weakly in the direction $u$. 
\item When $\N\colon\H\to\H$ is monotone (see definition \ref{def:monotony}(1)), for each $u\in\H$ and $t>s$ we have
\[
(t-s)\langle\N(tu)-\N(su),u\rangle_r\geqslant 0
\]
and thus the function $t\mapsto\langle\N(tu),u\rangle_r$ is nondecreasing. In particular the limit $\lim_{t\to+\infty}\langle\N(tu),u\rangle_r$ exists (perhaps it is $+\infty$).
\end{enumerate}
\end{remark}

\begin{lemma}\label{lem:recession}
Let $\N\colon\H\to\H$. Then
\begin{enumerate}[(1)]
\item $J_{\N}(\lam u)=\lam J_{\N}(u)$, for all $\lam>0$,
\item $J_{\N}(0)\in\{-\infty, 0\}$.
\end{enumerate}
Moreover, if $\N$ is monotone and $\psi_{\N}(u)\to+\infty$ as $\|u\|\to+\infty$, the functional $J_{\N}$ is weakly sequentially lower semicontinuous.
\end{lemma}

\begin{proof}
To prove \textit{(1)} note that for each $t_k\in\bbR$ and $v_k\in\H$
\[
\langle\N(t_kv_k),v_k\rangle=\lam\left\langle\N\left(\lam t_k\frac{1}{\lam} v_k\right),\frac{1}{\lam}v_k\right\rangle.
\]
Moreover $t_k\to+\infty$ and $v_k\weak\lam u$ if and only if $\lam t_k\to+\infty$ and $\frac{1}{\lam}v_k\weak u$.

Ad. \textit{(2)}. Taking $v_k=0$ for each $k\in\bbR$ we see that $J_{\N}(0)\leqslant 0$. From \textit{(1)} $J_{\N}(0)=\lam J_{\N}(0)$ for every $\lam>0$. Hence $J_{\N}(0)\in\{-\infty,0\}$.\smallskip

To prove the last part we will show that for each $c\in[-\infty,+\infty)$ the set
\[
S_c=\{u\in\H:\ J_{\N}(u)\leqslant c\}
\]
is weakly sequentially closed. To this end let $\{u_k\}_{k\in\bbN}\subset S_c$ be such that $u_k\weak u\in\H$. Then, for every $k\in\bbN$, $J_{\N}(u_k)\leqslant c$. Let $\{e_i\}_{i\in\bbN}$ be an orthonormal basis of $\H$ and recall that the weak convergence of a sequence in $\H$ is equivalent to the fact that it is bounded and the scalar products of its elements with each $e_i$ converge to the products with the limit. Hence for each $k\in\bbN$, from the definition of $J_{\N}(u_k)$, we can choose $t_k>k$ and $v_k\in\H$ such that
\begin{equation}\label{eq:recession_1}
|\langle e_i,v_k\rangle-\langle e_i,u_k\rangle|<\frac{1}{k},
\end{equation}
for all $i\in\{1,\ldots,k\}$, and
\begin{equation}\label{eq:recession_2}
\langle\N(t_kv_k),v_k\rangle_r\leqslant c+\frac{1}{k}\quad (\leqslant -k\ \text{if}\ c=-\infty).
\end{equation}
The sequence $\{v_k\}_{k\in\bbN}$ is bounded. Indeed, since $\N$ is monotone the function $t\mapsto\langle\N(tv_k),v_k\rangle_r$ is nondecreasing (see remark \ref{rem:recession}$(c)$) and we have from (\ref{eq:recession_2})
\[
\psi_{\N}(v_k)=\langle\N(v_k),v_k\rangle_r\leqslant\langle\N(t_kv_k),v_k\rangle_r\leqslant c+\frac{1}{k}\quad (\leqslant -k\ \text{if}\ c=-\infty)
\]
for each $k\in\bbN$. If $\|v_k\|\to+\infty$ we would have a contradiction. Moreover, for each $i\in\bbN$ and $k\geqslant i$ we have from (\ref{eq:recession_1})
\[
|\langle e_i,v_k\rangle-\langle e_i,u\rangle|\leqslant \frac{1}{k}+|\langle e_i,u_k\rangle-\langle e_i,u\rangle|
\]
and the right hand side goes to zero as $k$ goes to infinity. Therefore $t_k\to+\infty$ and $v_k\weak u$ thus from (\ref{eq:recession_2}) we get $J_{\N}(u)\leqslant c$.
\end{proof}
\begin{remark}\label{rem:recession2}
Let $\N\colon\H\to\H$ satisfy the assumptions of theorem \ref{thm:main_ex_monot}. We have noticed in remark \ref{rem:main_ex_monot} that condition $(i)$ of theorem \ref{thm:main_ex_monot} implies that $\N$ is monotone. Moreover, estimate $(ii)$ implies clearly that $\psi_{\N}$ is coercive. Hence the assumptions of lemma \ref{lem:recession} are satisfied and $J_{\N}$ is weakly sequentially lower semicontinuous.

We will show that condition $(ii)$ from theorem \ref{thm:main_ex_monot} implies that
\begin{equation}\label{eq:recession_3}
J_{\N}(u)\geqslant\frac{\ga\|h\|}{\del^2\al-\ga}\|u\|
\end{equation}
for each $u\in\ker\LL\bsh\{0\}$. Indeed, let $\veps>0$ and take sequences $t_k\to+\infty$ and $v_k\weak u$ such that
\[
\liminf_{k\to+\infty}\langle N(t_kv_k),v_k\rangle_r<J_{\N}(u)+\veps.
\]
Passing to a subsequence we can assume that we have $\|v_k\|\geqslant\|u\|-\veps>0$ for all $k\in\bbN$ (since $u\neq 0$) and that the limit inferior on the left hand side is actually the limit. If we put $u_k=t_kv_k$ then $\|u_k\|=t_k\|v_k\|\geqslant t_k(\|u\|-\veps)\to+\infty$ as $k\to+\infty$ so condition $(ii)$ implies
\[
\frac{\ga\|h\|}{\del^2\al-\ga}<\limsup_{k\to+\infty}\frac{\langle\N(u_k),u_k\rangle_r}{\|u_k\|}\leqslant \limsup_{k\to+\infty}\frac{\langle\N(u_k),u_k\rangle_r}{t_k(\|u\|-\veps)}= \frac{1}{\|u\|-\veps}\lim_{k\to+\infty}\langle N(t_kv_k),v_k\rangle_r< \frac{J_{\N}(u)+\veps}{\|u\|-\veps}.
\]
Since $\veps>0$ was arbitrary we get (\ref{eq:recession_3}).

Let us also note that lemma \ref{lem:recession}(1) and condition $(iii)$ imply for all $\|u\|\neq 0$
\[
J_{\N}(u)=\|u\|J_{\N}\left(\frac{u}{\|u\|}\right)>\|u\|\left(\frac{\ga\|h\|}{\del^2\al-\ga}+\left\langle h,\frac{u}{\|u\|}\right\rangle_r\right)
\]
so
\begin{equation}\tag*{\textit{(iii)$^\prime$}}
J_{\N}(u)>\frac{\ga\|h\|}{\del^2\al-\ga}\|u\|+\langle h,u\rangle_r.
\end{equation}
We see that $(iii)'$ improves estimate (\ref{eq:recession_3}) (which is a consequence of $(ii)$) precisely on a cone $C_h=\{u\in\H:\ \langle h,u\rangle_r>0\}$. Conditions $(iii)$ and $(iii)'$ are clearly equivalent for $\|u\|\geqslant 1$, however the author does not know if condition $(iii)$ can be replaced by $(iii)'$ in theorem \ref{thm:main_ex_monot}. Let us also note that according to lemma \ref{lem:recession}(2) condition $(iii)$ implies also that $J_{\N}(u)$ is not continuous at $0$.
\end{remark}
\medskip

Now we will recall the basic definitions and results in the theory of maximal monotone operators on Hilbert spaces. This theory is usually formulated in real Hilbert spaces but it can be easily extended to complex ones. Hence we will give all definitions and theorems assuming, like earlier, that $\H$ is complex.

Let $A\colon\H\to\power{\H}$ be a multi-valued operator (multifunction), where $\power{\H}$ denotes the family of all subsets of $\H$. The \emph{domain} of operator $A$ is the set
\[
\dom(A)=\{u\in\H:\ A(u)\neq\emptyset\},
\]
and its \emph{image} is defined as follows 
\[
\ran(A)=\bigcup_{u\in\H}A(u).
\]
If $A,B\colon\H\to\power{\H}$ and $\al,\be\in\bbR$ for all $u\in\H$ we set
\[
(\al A+\be B)(u)=\{\al v+\be v':\ v\in A(u),\ v'\in B(u)\}
\]
and then $\dom(\al A+\be B)=\dom(A)\cap\dom(B)$. The \emph{inverse operator} $A^{-1}\colon\H\to\power{\H}$ is defined for all $v\in\H$ by formula
\[
A^{-1}(v)=\{u\in\H:\ v\in A(u)\}.
\]
If we identify $A$ with its graph $\graph(A)=\{(u,v)\in\H\oplus\H:\ v\in A(u)\}$ then $A^{-1}$ is the operator whose graph is symmetric with respect to the graph of $A$, i.e., $(u,v)\in\graph(A)$ iff $(v,u)\in\graph(A^{-1})$. Of course we have $\dom(A^{-1})=\ran(A)$.\medskip

\begin{remark}
If for all $u\in\dom(A)$ the set $A(u)$ consists precisely of one element (so it is a singleton) we call $A$ a single-valued operator. In this case we can attribute to $u$ the unique element of $A(u)$, call it $\tilde{A}(u)$, so that we have a map $\tilde{A}\colon\dom(A)\supset\H\to\H$. Henceforward we will not distinguish between $A$ and $\tilde{A}$. 
\end{remark}\medskip

\begin{definition}\label{def:monotony}
We say that operator $A\colon\H\to\power{\H}$ is
\begin{enumerate}[(1)]
\item \emph{monotone} if for all $u,u'\in\dom(A)$
\[
\langle A(u)-A(u'),u-u'\rangle_r\geqslant 0,
\]
by which we mean that for all $v\in A(u),v'\in A(u')$ we have
\[
\langle v-v',u-u'\rangle_r\geqslant 0,
\]
\item \emph{strongly monotone} if there is a constant $C>0$ such that for all $u,u'\in\dom(A)$ we have
\[
\langle A(u)-A(u'),u-u'\rangle_r\geqslant C\|u-u'\|^2.
\]
\end{enumerate}
\end{definition}


In infinite-dimensional spaces the simplest examples of monotone operators are among linear and non-negative. Hence, according to lemma \ref{lem:decomp} operator $\rL$ is monotone on $\rH$. 

Since the family of monotone operators is inductive with respect to inclusion of graphs the following definition seems natural.

\begin{definition}\label{def:max_mon}
Assume that an operator $A\colon\H\to\power{\H}$ is monotone. Then we say that $A$ is \emph{maximal monotone} if for all $u\in\dom(A)$ and $u_0,v_0\in \H$ the condition
\[
\langle A(u)-v_0,u-u_0\rangle_r\geqslant 0,
\]
which we understand as
\[
\langle v-v_0,u-u_0\rangle_r\geqslant 0,\quad \text{for all}\ v\in A(u),
\]
implies that
\[
u_0\in\dom(A)\ \text{and}\ v_0\in A(u_0).
\]
\end{definition}\medskip

Maximal monotonicity means precisely that a given monotone operator does not have non-trivial monotone extension. Indeed, if we had $\langle A(u)-v_0,u-u_0\rangle_r\geqslant 0$ on $\dom(A)$ but $u_0\notin\dom(A)$ then defining operator $\tilde{A}\colon\dom(A)\cup\{u_0\}\to\power{\H}$ as $\tilde{A}(u)=A(u)$ for $u\in\dom(A)$, $\tilde{A}(u_0)=\{v_0\}$, we would get the non-trivial monotone extension. If $u_0\in\dom(A)$ but $v_0\notin A(u_0)$, we can define the monotone extension $\tilde{A}\colon\dom(A)\to 2^{\H}$ by formula $\tilde{A}(u)=A(u)$ for $u\in\dom(A)\bsh\{u_0\}$ and $\tilde{A}(u_0)=A(u_0)\cup\{v_0\}$.


The characterisation given in next theorem is fundamental in the study and applications of maximal monotone operators.
\begin{theorem}\label{thm:char_max_mon}
Let us assume that operator $A\colon\H\to\power{\H}$ is monotone. The following conditions are equivalent
\begin{enumerate}[$(i)$]
\item $A$ is maximal monotone,
\item $\ran(A+I\lam)=\H$ for some $\lam>0$,
\item $\ran(A+I\lam)=\H$ for all $\lam>0$.
\end{enumerate}
\end{theorem}
\begin{proof}
For the proof see for example \cite[Proposition 2.2, p. 23]{Brezis1973}.
\end{proof} 
\medskip

We will finish this subsection with some basic facts concerning maximal monotone operators. Firstly, we recall how this notion behaves under various operations on mappings.
\begin{lemma}\label{lem:pers_max_mon}
Assume that $A,B\colon\H\to\power{\H}$ are maximal monotone. Then
\begin{enumerate}[$(i)$]
\item for every $\lam>0$ operator $\lam A$ is maximal monotone,
\item the operator $A^{-1}$ is maximal monotone,
\item if \ $\textnormal{int}\dom(A)\cap\dom(B)\neq\emptyset$, operator $A+B$ is maximal monotone and
\[
\cl{\dom(A)\cap\dom(B)}=\cl{\dom(A)}\cap\cl{\dom(B)}.
\]
\end{enumerate}
\end{lemma}
\begin{proof}
Items $(i)$ and $(ii)$ follow easily from the definition of maximal monotonicity. For the proof of $(iii)$ see \cite[Cor. 2.7, p. 36]{Brezis1973}.
\end{proof}

Secondly, let us mention an easy result concerning surjectivity of maximal monotone operators (for more see \cite[p. 30-34]{Brezis1973}).
\begin{lemma}\label{lem:sur_max_mon}
If $A\colon\H\to\power{\H}$ is maximal and strongly monotone operator, then $\ran(A)=\H$.
\end{lemma}
\begin{proof}
Strong monotonicity (see definition \ref{def:monotony}$(2)$) implies that $A-C$ is also maximal monotone for some constant $C>0$. So conclusion follows from theorem \ref{thm:char_max_mon}.
\end{proof}

Finally, the following lemma binds together the theory of maximal monotone operators just sketched with the spectral assumptions imposed on linear operator $\LL$.
\begin{lemma}\label{lem:lin_max_mon}
Let $\rP=I-E_{0-}$ be the projection on $\rH$. Then operator $\rL\rP\colon\dom(\LL)\subset\H\to\rH$ is maximal monotone.
\end{lemma}
\begin{proof}
Since for any $\lam\in\bbC$ the operator $\rL\rP-\lam I$ is bijection from $\dom(\LL)$ onto $\H$ if and only if the operator $\rL-\lam I_+$, where $I_+$ is the identity on $\rH$, is bijection from $\dom(\rL)$ onto $\rH$ we have from lemma \ref{lem:decomp}
\[
\sp(\rL\rP)=\sp(\rL)\subset[0,+\infty).
\]
This shows that $\rL\rP$ is monotone and that for all $\lam>0$ we have $\ran(\rL\rP+\lam I)=\H$. Hence the conclusion follows from theorem \ref{thm:char_max_mon}.
\end{proof}

\section{Perturbed equations}

Firstly, we will show the solvability of the perturbed equations
\begin{equation}\label{eq:perturbed}
\veps\rP u+\LL u+\N(u)=h.
\end{equation}
From now on we will use the decomposition $\N(u)=\lN(u)+\rN(u)$ of the values of nonlinear part $\N$ according to the decomposition of Hilbert space $\H$ given in (\ref{eq:decomp}). For $u\in\H$ we will also write $u=u^-+u^+$ where $u^\pm\in\H_{\pm}$.


\begin{theorem}\label{thm:pert}
Assume that a self-adjoint operator $\LL\in\lin(\H)$ satisfies \ref{h:spectr}, \ref{h:gap} and \ref{h:bound}. Operator $\N\colon\H\to\H$ is bounded, demicontinuous, $\N(0)=0$ and
\begin{enumerate}[(i)]
\item there exists $\al>\ga/\del^2$ such that for all $u_1,u_2\in\H$ 
\[
\langle\N(u_1)-\N(u_2),u_1-u_2\rangle_r\geqslant\al\|\lN(u_1)-\lN(u_2)\|^2.
\]
\end{enumerate}
Then for each $\veps>0$ and $h\in\H$ the perturbed equation (\ref{eq:perturbed})
admits precisely one solution.
\end{theorem}

\begin{proof}
Let us denote by $\lP=E_{-\del}$ and $\rP=I-E_{0-}$ the orthoprojections on $\lH$ and $\rH$ respectively. Put
\[
\LL_\veps=\lL\lP+\veps\rP\in\bound(\H).
\]
The boundedness of $\lL$ was noted in lemma \ref{lem:decomp}. It is easy to check that $\LL_\veps$ is a symmetric bijection and hence $0\notin\sp(\LL_\veps)$. What is more, from the fact that $\sp(\lL)\subset[-\ga,-\del]$ (lemma \ref{lem:decomp}) we get
\[
\sp(\LL_\veps)\subset[-\ga,-\del]\cup\{\veps\}
\]
so $\LL_\veps^{-1}\in\bound(\H)$ is self-adjoint (because it is symmetric and has real spectrum) and applying spectral mapping theorem (see \cite[Prop. 5.5.3, p. 178]{BlaExnHav2008})
\[
\sp(\LL_\veps^{-1})\subset[-1/\del,-1/\ga]\cup\{1/\veps\}.
\]
Moreover if $u=\LL^{-1}_\veps v$ then $v^-=\lL u^- $ and $v^+=\veps u^+$ so we have
\begin{multline}\label{eq:thm_pert1}
\langle\LL^{-1}_\veps v, v\rangle=\langle u,\LL_\veps u\rangle=\langle u^-,\lL u^-\rangle+\veps\|u^+\|^2\geqslant\\
\geqslant1/\veps\|\veps u^+\|^2-\ga/\del^2\|\lL u^-\|^2=1/\veps\|v^+\|^2-\ga/\del^2\|v^-\|^2,
\end{multline}
where in the estimate from below we used the estimate in lemma \ref{lem:decomp}. Note that using the introduced notation we can convert equation (\ref{eq:perturbed}) as follows. Firstly we split the linear part
\[
\veps u^+ +\lL u^- = h-\N(u)-\rL u^+,
\]
by definition of $\LL_\veps$ we get
\[
\LL_\veps u=h-\N(u)-\rL u^+,
\]
we invert
\[
u=\LL^{-1}_\veps(h-\N(u)-\rL u^+),
\]
and finally get
\[
u=\LL^{-1}_\veps h-\LL^{-1}_\veps(\rL\rP+\N)(u).
\]
Let $v=(\rL\rP+\N)(u)$ and put
\[
A=(\rL\rP+\N)^{-1}\colon\H\supset\dom(A)\to 2^\H.
\]
According to the transformations just performed equation (\ref{eq:perturbed}) is equivalent to
\begin{equation}\label{eq:thm_pert2}
\LL^{-1}_\veps h\in A(v)+\LL^{-1}_\veps v,\quad v\in\dom(A).
\end{equation}

We showed in lemma \ref{lem:lin_max_mon} that $\rL\rP$ is maximal monotone. Because of $(i)$ operator $\N$ is in particular monotone. It is also maximal monotone. Indeed, assume that $u_0,v_0\in\H$ satisfy
\[
\langle N(u)-v_0,u-u_0\rangle_r\geqslant 0
\]
for every $u\in\H$. Then for all $u\in\H$ and $t>0$ we have in particular $\langle N(u_0+tu)-v_0,u\rangle_r\geqslant 0$ and making $t$ converge to $0$ we arrive at
\[
\langle N(u_0)-v_0,u\rangle_r\geqslant 0,
\]
since $N$ is demicontinuous. This gives $N(u_0)=v_0$ because $u$ is arbitrary. Hence operator $\N+\rL\rP$ is also maximal monotone (lemma \ref{lem:pers_max_mon}$(iii)$) and in consequence $A$ is maximal monotone (lemma \ref{lem:pers_max_mon}$(ii)$).

Let $v_1,v_2\in\dom(A)$. Then for all $u_1\in A(v_1),\ u_2\in A(v_2)$
\begin{multline*}
\langle u_1-u_2,v_1-v_2\rangle_r= \langle u_1-u_2,\rL u^+_1+\N(u_1)-\rL u^+_2-\N(u_2)\rangle_r= \langle u^+_1-u^+_2,\rL(u^+_1-u^+_2)\rangle+\\
+\langle u_1-u_2,\N(u_1)-\N(u_2)\rangle_r\geqslant\al\|\lN(u_1)-\lN(u_2)\|^2=\al\|v^-_1-v^-_2\|^2.
\end{multline*}
In the estimate from below we used the non-negativity of operator $\rL$ (see lemma \ref{lem:decomp}) and assumption $(i)$. Hence for all $v_1,v_2\in\dom(A)$
\begin{equation}\label{eq:thm_pert3}
 \langle A(v_1)-A(v_2),v_1-v_2\rangle_r\geqslant\al\|v^-_1-v^-_2\|^2.
\end{equation}
Put $A_\veps=A+\LL^{-1}_\veps$. From (\ref{eq:thm_pert1}) and (\ref{eq:thm_pert3}) we get for each $v_1,v_2\in\dom(A)$ that
\begin{multline*}
 \langle A_\veps(v_1)-A_\veps(v_2),v_1-v_2\rangle_r= \langle A(v_1)-A(v_2),v_1-v_2\rangle_r+\langle \LL^{-1}_\veps(v_1-v_2),v_1-v_2\rangle\geqslant\\
\geqslant\al\|v^-_1-v^-_2\|^2+1/\veps\|v^+_1-v^+_2\|^2-\ga/\del^2\|v^-_1-v^-_2\|^2= 1/\veps\|v^+_1-v^+_2\|^2+(\al-\ga/\del^2)\|v^-_1-v^-_2\|^2.
\end{multline*}
So operator $A_\veps$ is strongly monotone (definition \ref{def:monotony}$(2)$) with constant $C=\min\{1/\veps,\al-\ga/\del^2\}>0$. In particular it means that it is one-to-one in a sense that if $v_1\neq v_2$ then $A_\veps(v_1)\cap A_\veps(v_2)=\emptyset$. In fact, if there was a vector $w\in A_\veps(v_1)\cap A_\veps(v_2)$ then from the last estimate we would get
\[
0=\langle w-w,v_1-v_2\rangle\geqslant C\|v_1-v_2\|^2.
\]

Now we will show that operator $A_\veps$ is maximal monotone. To this end note that
\[
A_\veps+\lam=A+(\LL^{-1}_\veps+\lam)
\]
is surjective for all $\lam>1/\del$. Indeed, operator $A$ is, as we already proved, maximal monotone and $\LL_\veps^{-1}+1/\del$ also belongs to this class according to theorem \ref{thm:char_max_mon} and the fact that $(-\infty,-1/\del)\subset\rez(\LL^{-1}_\veps)$. From lemma \ref{lem:pers_max_mon}$(iii)$ we get that $A+\LL^{-1}_{\veps}+1/\delta$ is maximal monotone so, making use of theorem \ref{thm:char_max_mon} once more, for every $\lambda>1/\delta$ the operator $A+\LL^{-1}_{\veps}+\lam$ is surjective. 

To sum up, we showed that operator $A_\veps$ is strongly and maximal monotone so according to lemma \ref{lem:sur_max_mon} it is surjective. This ensures the existence of solution to (\ref{eq:thm_pert2}) which is unique since set-values of $A_\veps$ are disjoint, this ends the proof of the theorem.
\end{proof}

\section{Proof of theorem \ref{thm:main_ex_monot}}

Firstly we will show that the boundedness of solution of perturbed equations, when the perturbation parameter $\veps$ is also bounded, is sufficient for the solvability of the main equation.

\begin{lemma}\label{lem:limit_sol_monot}
Assume that $\LL$ and $\N$ satisfy the assumptions of theorem \ref{thm:pert} and let $u_\veps\in\H$ be a solution of perturbed equation (\ref{eq:perturbed}) with $\veps>0$. If there is a constant $C>0$ such that $\|u_\veps\|\leqslant C$ for each sufficiently small $\veps>0$ then equation (\ref{eq:main}) has a solution.
\end{lemma}

\begin{proof}
Fix a sequence $\veps_k\to 0$ and let $u_k:=u_{\veps_k}\in\H$ be a solution of equation \ref{eq:perturbed} with $\veps=\veps_k$. We can assume, choosing a subsequence if needed, that there is $u\in\H$ such that $u_k\weak u$. Since $\N$ is monotone we have for all $v\in\H$
\[
\langle \N(u_k)-\N(v),u_k-v\rangle_r\geqslant 0,
\]
and further
\[
 \langle h-\veps_ku^+_k-\LL u_k-\N(v),u_k-v\rangle_r\geqslant 0.
\]
Hence, taking advantage of non-negativity of $\rL$, we have for all $v\in\dom(\LL)$
\begin{multline}\label{eq:lem_pert}
\langle h-\N(v),u_k-v\rangle_r\geqslant
 \langle\LL u_k+\veps_ku^+_{k},u_k-v\rangle_r= \langle\lL u^-_{k},u^-_{k}-v^-\rangle_r+\\
+\langle\rL(u^+_{k}-v^+),u^+_{k}-v^+\rangle + 
\langle\rL v^+,u^+_{k}-v^+\rangle_r+ \langle\veps_ku^+_{k},u^+_{k}-v^+\rangle_r\geqslant\\
\geqslant\langle\lL u^-_{k},u^-_{k}-v^-\rangle_r+ 
\langle\rL v^+,u^+_{k}-v^+\rangle_r + \langle\veps_ku^+_{k},u^+_{k}-v^+\rangle_r.
\end{multline}
Since $\veps_ku^+_{k}\to 0$ and $u^+_{k}\weak u^+$ it follows that
\[
\lim_{k\to\infty}\langle\veps_ku^+_{k},u^+_{k}-v^+\rangle\to 0
\]
and
\[
\lim_{k\to\infty}\langle\rL v^+,u^+_{k}-v^+\rangle= 
\langle\rL v^+,u^+-v^+\rangle.
\]
We have to study the convergence of the first term on the right hand side of (\ref{eq:lem_pert}). Firstly, we will show that the sequence $\{\lL u^-_{k}\}_{k\in\bbN}$ is convergent in norm. To this end fix $k,l\in\bbN$ and note that making use of assumption $(i)$ from theorem \ref{thm:pert} we have
\begin{multline*}
\|\lL u^-_{k}-\lL u^-_{l}\|^2=\|\lN(u_k)-\lN(u_l)\|^2\leqslant 1/\al\langle \N(u_k)-\N(u_l),u_k-u_l\rangle_r=\\
=-1/\al \langle\LL u_k-\LL u_l,u_k-u_l\rangle+ 1/\al \langle\veps_lu^+_{l}-\veps_ku^+_{k},u^+_{k}-u^+_{l}\rangle_r.
\end{multline*}
The non-negativity of $\rL$, lemma \ref{lem:decomp}$(i)$ and condition \ref{h:bound} imply that
\begin{multline*}
 \langle\LL u_k-\LL u_l,u_k-u_l\rangle= \langle\lL(u^-_{k}-u^-_{l}),u^-_{k}-u^-_{l}\rangle+ \langle\rL(u^+_{k}-u^+_{l}),u^+_{k}-u^+_{l}\rangle\geqslant\\
\geqslant-\ga/\del^2\|\lL(u^-_{k}-u^-_{l})\|^2.
\end{multline*}
From this two estimates we eventually get
\[
\|\lL u^-_{k}-\lL u^-_{l}\|^2\leqslant\ga/\al\del^2\|\lL(u^-_{k}-u^-_{l})\|^2 + C(\veps_k+\veps_l),
\]
and since $\al>\ga/\del^2$ we compute
\[
\lim_{k,l\to\infty}(1-\ga/\al\del^2)\|\lL u^-_{k}-\lL u^-_{l}\|^2\leqslant \lim_{k,l\to\infty}C(\veps_k+\veps_l)=0.
\]
We know from lemma \ref{lem:decomp} that $\lL\in\bound(\lH)$ so $\lL u^-_{k}\weak \lL u^-$, according to the weak convergence of $\{u^-_{k}\}_{k\in\bbN}$, and since $\{\lL u^-_{k}\}_{k\in\bbN}$ is a Cauchy sequence in norm it follows that $\lL u^-_{k}\to\lL u^-$. Therefore we have
\[
\lim_{k\to\infty}\langle\lL u^-_{k},u^-_{k}-v^-\rangle= 
\langle\lL u^-,u^--v^-\rangle.
\]
Taking the limit on the both sides in (\ref{eq:lem_pert}) as $k\to\infty$ we arrive at
\begin{equation}\label{eq:lem_pert2}
\langle h-\N(v),u-v\rangle_r\geqslant\langle\lL u^-,u^--v^-\rangle_r + \langle\rL v^+,u^+-v^+\rangle_r,
\end{equation} 
for every $v\in\dom(\LL)$.

Putting $v^-=u^-$ in (\ref{eq:lem_pert2}) we have for all $v^+\in\dom(\rL)$
\[
 \langle\rL v^++\N(u^-+v^+)-h^+,v^+-u^+\rangle_r\geqslant 0
\]
which, from the maximal monotonicity of $\rL(\cdot)+\N(u^-+\cdot\,)$, gives $u^+\in\dom(\rL)$ and
\[
\rL u^++\N(u^-+u^+)=h^+.
\]

Next, fixing $v^+=u^+$ we get from (\ref{eq:lem_pert2})
\[
\langle h^--\lL u^--\lN(v^-+u^+),u^--v^-\rangle_r\geqslant 0.
\]
If we take $v^-=u^-+tw^-\in\lH$, where $w^-\in\lH$ and $t>0$, we then have
\[
\langle h^--\lL u^--\lN(u+tw^-),w^-\rangle_r\geqslant 0.
\]
Making $t$ converge to zero and using demi-continuity of $\N$ it follows that
\[
\langle h^--\lL u^--\lN(u),w^-\rangle_r\geqslant 0
\]
which, because $w^-\in\lH$ is arbitrary, implies $\lL u^-+\lN(u)=h^-$.
\end{proof}

We are in position to prove theorem \ref{thm:main_ex_monot}. Recall that by $J_{\N}$ we denote the recession functional of operator $\N$ with respect to weak convergence introduced in definition \ref{def:recession}.

\begin{proof}[Proof of theorem \ref{thm:main_ex_monot}]
Let $\veps_k\to 0$ and let $u_k\in\dom(\LL)$ be a solution of perturbed equation (\ref{eq:perturbed}) with $\veps=\veps_k$. According to lemma \ref{lem:limit_sol_monot} it suffices to show that the sequence $\{u_k\}_{k\in\bbN}$ is bounded.\medskip

\noindent Step 1. Firstly, we will show that from boundedness of sequence $\{u^+_{k}\}_{k\in\bbN}$ follows the boundedness of $\{u^-_{k}\}_{k\in\bbN}$. Therefore assume that the sequence $\{u^+_{k}\}_{k\in\bbN}$ is bounded. Since $u_k,\ k\in\bbN$, is a solution of perturbed equation we have in particular
\[
u^-_{k}=K(h^--\lN(u_k)),
\]
where $K=\lL^{-1}$ belongs to $\bound(\lH)$ (see lemma \ref{lem:decomp}). Thus it is sufficient to show that the sequence $\{\lN(u_k)\}_{k\in\bbN}$ is bounded. We will prove more, i.e, that the sequence $\{\N(u_k)\}_{k\in\bbN}$ is bounded. To this end let us note that multiplying both sides of (\ref{eq:perturbed}) with $u_k$ we get
\[
\langle\veps_ku^+_{k}+\LL u_k+\N(u_k),u_k\rangle=\langle h,u_k\rangle,
\]
hence
\[
\langle\N(u_k)-h,u_k\rangle=-\langle\LL u_k,u_k\rangle-\veps_k\|u^+_{k}\|^2.
\]
On the other hand, making use of assumption $(i)$ with $u_1=u_k$ and $u_2=0$ we have
\[
\langle\N(u_k)-h,u_k\rangle=\langle\N(u_k),u_k\rangle_r-\langle h,u_k\rangle_r \geqslant \al\|\N(u_k)\|^2-\|h\|\|u_k\|.
\]
Taking together the two above formulas we arrive at
\[
\al\|\N(u_k)\|^2-\|h\|\|u_k\|\leqslant-\langle\LL u_k,u_k\rangle-\veps_k\|u^+_{k}\|^2 \leqslant\ga/\del^2\|\lL u^-_{k}\|^2-\veps_k\|u^+_{k}\|^2,
\]
where, in the last estimate, we used lemma \ref{lem:decomp}. Next using in the first place $\|u^-_{k}\|\leqslant1/\del\|\lL u^-_{k}\|$ and afterwards $\|\lL u^-_{k}\|\leqslant\|\N(u_k)\|+\|h\|$ we compute
\begin{multline*}
\al\|\N(u_k)\|^2\leqslant\ga/\del^2\|\lL u^-_{k}\|^2+\|h\|\|u^-_{k}\|+\|h\|\|u^+_{k}\|-\veps_k\|u^+_{k}\|^2\leqslant \ga/\del^2\|\lL u^-_{k}\|^2+\\[1.5ex]
+1/\del\|h\|\|\lL u^-_{k}\|+\|h\|\|u^+_{k}\|-\veps_k\|u^+_{k}\|^2=\ga/\del^2\big(\|\lL u^-_{k}\|+\del\|h\|/(2\ga)\big)^2+\|h\|\|u^+_{k}\|-\veps_k\|u^+_{k}\|^2-\\[1.5ex]
-\|h\|^2/(2\ga)^2\leqslant\ga/\del^2\big(\|\N(u_k)\|+\|h\|+\del\|h\|/(2\ga)\big)^2+\|h\|\|u^+_{k}\|-\veps_k\|u^+_{k}\|^2-\|h\|^2/(2\ga)^2.
\end{multline*}
Taking advantage of Cauchy inequality with $\veps'>0$ such that $(1+2\veps')\ga/\del^2<\al$ we get
\[
\al\|\N(u_k)\|^2\leqslant(1+2\veps')\ga/\del^2\|\N(u_k)\|^2+\|h\|\|u^+_{k}\|-\veps_k\|u^+_{k}\|^2+C(\veps'),
\]
where $C(\veps')$ depends also on $\del$, $\ga$ and $\|h\|$, and  we finally arrive at
\begin{equation}\label{eq:thm_main_ex_monot1}
(\al-(1+2\veps')\ga/\del^2)\|\N(u_k)\|^2\leqslant\|h\|\|u^+_{k}\|-\veps_k\|u^+_{k}\|^2+C(\veps')
\end{equation}
from which the boundedness of $\{\N(u_k)\}_{k\in\bbN}$ follows.\medskip

\noindent Step 2. Therefore we suppose now that $\|u^+_{k}\|\to+\infty$ and search for contradiction. Firstly we will show that in this case
\begin{equation}\label{eq:thm_main_ex_monot2}
\limsup_{k\to\infty}\frac{\|u^-_{k}\|^2}{\|u^+_{k}\|}\leqslant\frac{\|h\|}{\del^2\al-\ga}.
\end{equation}
Indeed, since $u^-_{k}=K(h^--\lN(u_k))$ Cauchy inequality with $\veps'$ and estimate (\ref{eq:thm_main_ex_monot1}) imply that
\begin{multline*}
\frac{\|u^-_{k}\|^2}{\|u^+_{k}\|}\leqslant\frac{(1+2\veps')\|\N(u_k)\|^2+C(\veps')\|h\|^2}{\del^2\|u^+_{k}\|}\leqslant\\ \leqslant\frac{(1+2\veps')\|h\|}{\del^2\al-(1+2\veps')\ga}-\frac{(1+2\veps')\veps_k}{\del^2\al-(1+2\veps')\ga}\|u^+_{k}\|+\frac{C(\veps')(1+\|h\|^2)}{\|u^+_{k}\|},
\end{multline*}
and we get
\[
\frac{\|u^-_{k}\|^2}{\|u^+_{k}\|}\leqslant\frac{(1+2\veps')\|h\|}{\del^2\al-(1+2\veps')\ga}+\frac{C(\veps')(1+\|h\|^2)}{\|u^+_{k}\|}.
\]
Making $k$ converge to $\infty$ and then $\veps'$ to $0$ we arrive at (\ref{eq:thm_main_ex_monot2}). In particular (\ref{eq:thm_main_ex_monot2}) implies that
\[
\lim_{k\to\infty}\frac{\|u^-_{k}\|}{\|u^+_{k}\|}=0.
\]
Thus we can assume that there is $u\in\rH$ such that 
\[
v_k:=\frac{u_{k}}{\|u^+_{k}\|}\weak u.
\]


\noindent Step 3. Note that from \ref{h:bound} and (\ref{eq:thm_main_ex_monot2}) we have
\begin{equation*}
\limsup_{k\to\infty}\frac{\langle\lN(u_k),u^-_{k}\rangle_r}{\|u^+_{k}\|}=\limsup_{k\to\infty}\frac{\langle h^--\lL u^-_{k},u^-_{k}\rangle_r}{\|u^+_{k}\|}\leqslant \limsup_{k\to\infty}\frac{\langle h^-,u^-_{k}\rangle_r +\ga\|u^-_{k}\|^2}{\|u^+_{k}\|}\leqslant\frac{\ga\|h\|}{\del^2\al-\ga}.
\end{equation*}
Since operator $\LL+\veps_kI$ is non-negative on $\rH$ we also have
\begin{equation*}
\limsup_{k\to\infty}\frac{\langle\rN(u_k),u^+_{k}\rangle_r}{\|u^+_{k}\|}\leqslant\limsup_{k\to\infty}\frac{\langle(\rL+\veps_k)u^+_{k}+\rN(u_k),u^+_{k}\rangle_r}{\|u^+_{k}\|}=\limsup_{k\to\infty}\frac{\langle h^+,u^+_{k}\rangle_r}{\|u^+_{k}\|}=\langle h,u\rangle_r.
\end{equation*}
From these two estimates it follows that
\begin{equation}\label{eq:thm_main_ex_monot3}
\limsup_{k\to\infty}\frac{\langle\N(u_k),u_{k}\rangle_r}{\|u^+_{k}\|}\leqslant \frac{\ga\|h\|}{\del^2\al-\ga}+\langle h,u\rangle_r.
\end{equation}
Let us consider two cases. If $u=0$ then from (\ref{eq:thm_main_ex_monot3})
\[
\limsup_{k\to\infty}\frac{\langle\N(u_k),u_{k}\rangle_r}{\|u_{k}\|}\leqslant \frac{\ga\|h\|}{\del^2\al-\ga} 
\]
which is in contradiction with $(ii)$. If $u\neq 0$ then using (\ref{eq:thm_main_ex_monot3}) once more we get
\[
\frac{\ga\|h\|}{\del^2\al-\ga}+\langle h,u\rangle_r\geqslant  \liminf_{k\to\infty}\left\langle\N\left(\|u^+_{k}\|v_k\right),v_k\right\rangle_r\geqslant J_{\N}(u),
\]
which contradicts $(iii)$ this time.
\end{proof}

\section{Proof of theorem \ref{thm:Schrodinger_ex}}

Let $f\colon\bbR[n]\times\bbR\to\bbR$ be a Carath\'{e}odory function and for $u\in\LR[n]$ define
\begin{equation}\label{eq:superpos_op}
N(u)(x)=f(x,u(x)).
\end{equation}
Throughout this section $\langle\cdot\,,\cdot\rangle_{2}$ and $\|\cdot\|_{2}$ will denote respectively the scalar product and norm in $\LC[n]$, i.e.
\begin{align*}
&\langle u,v\rangle_{2}=\int_{\bbR[n]}u(x)\overline{v(x)}\,dx\\
&\|u\|_{2}=\int_{\bbR[n]}|u(x)|^2\,dx
\end{align*}
for all $u,v\in\LC[n]$. We treat $\LR[n]$ as a subspace of $\LC[n]$.

\begin{proof}[Proof of theorem \ref{thm:Schrodinger_ex}]
The estimate from above in assumption $(ii)$ ensures that the superposition operator $\N$, defined in (\ref{eq:superpos_op}), acts in $\LR[n]$, is bounded and continuous. It is even necessary for this to happen in the class of Carath\'{e}odory functions (or more generally for sup-measurable functions, see \cite[Theorem 3.1, p. 67]{AppellZabrejko1990}).\smallskip

\noindent Step 1. Firstly, we will show that for every $u,u'\in\LR[n]$ we have
\begin{equation*}
\langle\N(u)-N(u'),u-u'\rangle_{2}\geqslant\al\|\N(u)-\N(u')\|_{2}^2
\end{equation*}
and 
\[
\langle\N(u),u\rangle_{2}\geqslant a\|u\|_{2}^2.
\]
From $(iii)$ we know that for a.e. $x\in\bbR[n]$ and any $t,t'\in\bbR$ we have
\begin{equation}\label{eq:ex_monot1}
|f(x,t)-f(x,t')|\leqslant 1/\al|t-t'|.
\end{equation}
Multiplying both sides by $|f(x,t)-f(x,t')|$ and using the monotonicity of $f(x,\cdot\,)$ we get
\begin{equation}\label{eq:ex_monot2}
\al|f(x,t)-f(x,t')|^2\leqslant (f(x,t)-f(x,t'))(t-t').
\end{equation}
In an analogous fashion, using the fact that $f(x,t)t\geqslant 0$ we get from lower estimate in $(ii)$ that
\begin{equation}\label{eq:ex_monot3}
f(x,t)t\geqslant at^2.
\end{equation}
Finally, let us take $u,u'\in\LR[n]$. Firstly, applying (\ref{eq:ex_monot2}) we compute 
\begin{multline*}
\langle\N(u)-N(u'),u-u'\rangle_{2}=\int_{\bbR[n]}(f(x,u(x))-f(x,u'(x))(u(x)-u'(x))\,dx\geqslant\\
\geqslant\int_{\bbR[n]}\al|f(x,u(x))-f(x,u'(x))|^2\,dx=\al\|\N(u)-\N(u')\|_{2}^2,
\end{multline*}
which gives the former estimate. Secondly, (\ref{eq:ex_monot3}) implies that
\[
\langle\N(u),u\rangle_{2}=\int_{\bbR[n]}f(x,u(x))u(x)\,dx\geqslant a\|u\|_{2}^2,
\]
so the latter one is also true.\smallskip

\noindent Step 2. Condition (\ref{eq:Kato}) ensures that the operator $S=-\Del+V$ is well defined on $\dom(S)=H^{2}(\bbR[n],\bbC)$ and bounded from below. Moreover conditions \ref{h:spectr} -- \ref{h:bound} imposed on linear operator in theorem \ref{thm:main_ex_monot} are satisfied. Now let us define operator $\tilde{\N}\colon\LC[n]\to\LC[n]$ as follows
\[
\tilde{\N}(\xi)(x)=f(x,u(x))-if(x,v(x))=\N(u)(x)-i\N(v)(x),
\]
where $u,v\in\LR[n]$ and $\xi=u+iv\in\LC[n]$. Then we have
\[
\langle\tilde{\N}(\xi)-\tilde{\N}(\xi'),\xi-\xi'\rangle_{2,r}=\langle\N(u)-N(u'),u-u'\rangle_{2}+\langle\N(v)-N(v'),v-v'\rangle_{2}
\]
and
\[
\langle\tilde{\N}(\xi),\xi\rangle_{2,r}=\langle\N(u),u\rangle_{2}+ \langle\N(v),v\rangle_{2}
\]
for all $\xi,\xi'\in\LC[n]$, where $\langle\cdot\,,\cdot\rangle_{2,r}=\re\langle\cdot\,,\cdot\rangle_{2}$. Therefore from step 1 follows that $\tilde{\N}$ is continuous, bounded and
\[
\langle\tilde{\N}(\xi)-\tilde{\N}(\xi'),\xi-\xi'\rangle_{2,r}\geqslant\al \left(\|\N(u)-\N(u')\|_{2}^2+\|\N(v)-\N(v')\|_{2}^2\right)= \al\|\tilde{\N}(\xi)-\tilde{\N}(\xi')\|_{2}^2,
\]
and
\[
\limsup_{\|\xi\|_{2}\to\infty}\frac{\langle\tilde{\N}(\xi),\xi\rangle_{2,r}}{\|\xi\|_{2}}\geqslant\limsup_{\|\xi\|_{2}\to\infty}a\|\xi\|_{2}=+\infty.
\]
Hence assumptions $(i)$ and $(ii)$ of theorem \ref{thm:main_ex_monot} are satisfied. Next take $0<t_k\to\infty$ and $\{\eta_k\}\subset\LC[n]$ such that $\eta_k\weak \xi\neq 0$. Then $\liminf\|\eta_k\|\geqslant\|\xi\|=1$ so $t_k\|\eta_k\|\to+\infty$. Making use of condition (\ref{eq:ex_monot3}) we get
\[
\langle\tilde{\N}(t_k\eta_k),\eta_k\rangle_{2,r}=\int_{\bbR[n]}f(x,t_ku_k(x))u_k(x)\,dx+\int_{\bbR[n]}f(x,t_kv_k(x))v_k(x)\,dx\geqslant at_k\|\eta_k\|^{2}_{2}\to\infty,
\]
where for each $k\in\bbN$ we have $\eta_k=u_k+iv_k$. Hence
\[
J_{\tilde{\N}}(\xi)=+\infty
\]
for each $\xi\neq 0$, so condition $(iii)$ from theorem \ref{thm:main_ex_monot} is also satisfied. This implies that there is $\xi=u+iv\in\dom(S)=H^2(\bbR[n]\!,\bbC)$ such that
\[
S\xi(x)+\tilde{\N}(\xi)(x)=h(x)
\]
and in particular, since $h\in\LR[n]$, the real part $u\in H^{2}(\bbR[n],\bbR)$ satisfies
\[
S u(x)+f(x,u(x))=h(x)
\]
in $\LR[n]$.
\end{proof}

\bibliographystyle{plain}
\bibliography{C:/Users/pc/Bibliografia/Articles,C:/Users/pc/Bibliografia/Functional_Analysis,C:/Users/pc/Bibliografia/Nonlinear_Analysis,C:/Users/pc/Bibliografia/Mathematical_Physics}

\end{document}